\documentclass[reqno]{amsart}

\usepackage[utf8, latin1]{inputenc}
\usepackage[english]{babel}
\usepackage{stmaryrd}
\usepackage{amsthm}
\usepackage{amsmath}
\usepackage{amssymb}
\usepackage{enumitem}
\usepackage[dvips]{graphicx}
\usepackage{color}
\usepackage{hyperref}
\usepackage{url}
\usepackage[left=3.5cm, right=3.5cm]{geometry}
\usepackage{fancyhdr}
\usepackage{mathrsfs}
\usepackage{bm}
\usepackage{stmaryrd}

\newtheorem{theorem}{Theorem}

\newtheorem{corollary}{Corollary}
\newtheorem{rmk}{Remark}
\let\oldrmk\rmk
\renewcommand{\rmk}{\oldrmk\normalfont}
\newtheorem{example}{Example}
\let\oldexample\example
\renewcommand{\example}{\oldexample\normalfont}

\theoremstyle{definition}

\theoremstyle{remark}
\newtheorem{claim}{\textsc{Claim}}

\hypersetup{
    pdftitle={Sums of Multivariate Polynomials in Finite Subgroups},
    pdfauthor={Paolo Leonetti e Andrea Marino},
    pdfmenubar=false,
    pdffitwindow=true,
    pdfstartview=FitH,
    colorlinks=true,
    linkcolor=red,
    citecolor=green,
    urlcolor=cyan
}

\renewcommand{\emptyset}{\varnothing}

\renewcommand{\rho}{\varrho}

\newcommand*{\bigchi}{\mbox{\Large$\chi$}}
\newcommand{\free}{\natural}
\newcommand{\exponent}{\lambda}
\newcommand{\partition}{\mathscr{P}{\scriptsize \hspace{-.3mm}\text{art}}}
\newcommand{\interval}{\mathscr{I}{\scriptsize \hspace{-.5mm}\text{nt}}}

\begin{document}
\title{Sums of Multivariate Polynomials in Finite Subgroups}    

\author{Paolo Leonetti}
\address{Universit\`a L. Bocconi, via Roentgen $1$, $20136$ Milano, Italy.}
\email{leonetti.paolo@gmail.com}

\author{Andrea Marino}
\address{Scuola Normale Superiore, Collegio Timpano, Lungarno Pacinotti $51$, $56126$ Pisa, Italy.}
\email{andreamarino95@gmail.com}



\begin{abstract}  
\noindent{} Let $R$ be a commutative ring, $f \in R[X_1,\ldots,X_k]$ a multivariate polynomial, and $G$ a finite subgroup of the group of units of $R$ satisfying a certain constraint, which always holds if $R$ is a field. Then, we evaluate $\sum f(x_1,\ldots,x_k)$, where the summation is taken over all pairwise distinct $x_1,\ldots,x_k \in G$. In particular, let $p^s$ be a power of an odd prime, $n$ a positive integer coprime with $p-1$, and $a_1,\ldots,a_k$ integers such that $\varphi(p^s)$ divides $a_1+\cdots+a_k$ and $p-1$ does not divide $\sum_{i \in I}a_i$ for all non-empty proper subsets $I\subseteq \{1,\ldots,k\}$; then
\begin{displaymath}
\sum x_1^{a_1}\cdots x_k^{a_k} \equiv \frac{\varphi(p^s)}{\mathrm{gcd}(n,\varphi(p^s))}(-1)^{k-1}(k-1)! \hspace{1mm}\bmod{p^s}, 
\end{displaymath}
where the summation is taken over all pairwise distinct $n$-th residues $x_1,\ldots,x_k$ modulo $p^s$ coprime with $p$. 
\end{abstract}

\maketitle

\thispagestyle{empty} 

\section{Introduction}\label{sec:intro}

The aim of this article is to evaluate certain types of symmetric sums of distinct elements taken in a subgroup of the group of units of a commutative ring (hereafter, rings are always assumed to be unital and non-trivial). This includes, for instance, the case of sums of multivariate polynomials taking distinct values in the set of $n$-th residues modulo a prime. The first result of this type was obtained by Pierce \cite{Pierce}, who proved that an integral symmetric homogeneous function of degree $d$ of the $n$-th residues of an odd prime $p$ is divisible by $p$ if $d$ is not divisible by $(p-1)/\text{gcd}(n,p-1)$. Here, in particular, we evaluate sums of polynomials of $n$-th residues in the remaining case.

Moreover, symmetric sums of functions taking distinct values in a given set have been already studied in the literature: Ferrers \cite{Ferrers} proved the folklore result that an odd prime $p$ divides the sum of the products of the numbers $1,\ldots,p-1$, taken $k$ together, whenever $k<p-1$. Afterwards, this theorem was increasingly generalized by Glaisher \cite{Glaisher}, Moritz \cite{Moritz}, and Ricci \cite{Ricci}. All these results provide, in turn, generalizations of Wilson's theorem. Within this context, the proof of the celebrated Erd\H{o}s--Ginzburg--Ziv theorem \cite{EGZ} provided by Gao \cite{Gao} shows a clear connection between zero-sum problems in additive number theory and the study of sums of symmetric functions. 

Here below, let $R$ be a commutative ring and $G$ a finite subgroup of invertible elements with $n:=|G|$ and denote by $\lambda$ the exponent of $G$. (We refer to \S{} \ref{sec:notation} for notations used, but not defined, here.) Given a multivariate polynomial $f \in R[X_1,\ldots,X_k]$ and $k \in \mathbf{N}^+$ with $k\le n$, we provide a method to evaluate symmetric sums of the form
\begin{displaymath}
{\sum_{\substack{x_1,\ldots,x_k \,\in\, G \\ x_1,\ldots,x_k \hspace{1mm}\mathrm{pairwise}\hspace{1mm}\mathrm{distinct}}}}f(x_1,\ldots,x_k),
\end{displaymath}
whenever a certain condition which relates the structure of $G$ and the regular elements of $R$ is satisfied.

Since $f$ can be written as a finite sum of monomials $rX_1^{a_1}\cdots X_k^{a_k}$ for some $r \in R$ and $a_1,\ldots,a_k \in \mathbf{N}$, it is enough to evaluate symmetric sums of the form 
\begin{equation}\label{eq:pA}
p(A):={\sum_{\substack{x_1,\ldots,x_k\,\in\,G \\ x_1,\ldots,x_k \hspace{1mm}\mathrm{pairwise}\hspace{1mm}\mathrm{distinct}}}}x_1^{a_1}\cdots x_k^{a_k}.
\end{equation}

Here, $A$ stands for the multiset of integer exponents $\{a_1,\ldots,a_k\}$ (note that the order of the elements of $A$ does not matter since \eqref{eq:pA} is symmetric in $x_1,\ldots,x_n$). 

It is worth noting that, by Lagrange's theorem (see e.g. \cite{Mollin}), the order of each elements $g\in G$ divides $|G|$, hence $\exponent$ divides $|G|$. However, if $R$ is the ring $\mathbf{Z}_m$ of integers modulo $m$ and $G$ its subgroup of units, then the exponent of $G$ and $\varphi(m)$ (that is $|G|$) have different normal orders, see Erd\H{o}s, Pomerance and Schmutz \cite{EPS}.

\section{Notation}\label{sec:notation}

Through the paper, $\mathbf{Z}$, $\mathbf{N}$, and $\mathbf{N}^+$ stand, respectively, for the set of integers, non-negative integers, and positive integers. We use $\mathbf{Z}_m$ for the ring of integers modulo $m$.

Unless stated otherwise, $R$ stands always for a (non-trivial) commutative ring with a (non-zero) multiplicative identity, denoted by $1$. In this respect, let $D$ be the set of \emph{non-regular} elements of $R$ (we recall that $r \in R$ is said to be non-regular if there exists a non-zero $t \in R$ such that $rt=0$; in particular, $0 \in D$).

Then, $G$ denotes a (non-empty) finite subgroup of the group of units of $R$. The \emph{order} of each $g \in G$ is $\mathrm{ord}(g):=\min \{n \in \mathbf{N}^+: g^n=1\}$, while the \emph{exponent of} $G$, denoted by $\lambda$, is the least common multiple of $\{\mathrm{ord}(g): g \in G\}$ (however, it is easy to see that, since $G$ is a finite abelian group, then there exists $g\in G$ such that $\lambda= \mathrm{ord}(g)$).

We assume by convention that empty sums are equal to $0$.
Given a finite non-empty multiset $X$ of integers (that is, a set where repetitions are allowed), we define the \emph{sum of its elements} by $s(X):=\sum_{x \in X}x$. Accordingly, given the multiset of integer exponents $A$ and a subset $B\subseteq A$, let $\mathscr{P}(B)$ represent the collection of the partitions $\mathcal{P}$ of $B$ such that $\exponent$ divides $s(P)$ for all $P \in \mathcal{P}$ (in particular, $\mathscr{P}(B)=\emptyset$ if $\exponent$ does not divide $s(B)$). Lastly, define the \emph{characteristic number} of $B$ by
\begin{displaymath}
\bigchi(B):=|G|(-1)^{|B|-1}(|B|-1)!. 
\end{displaymath}

We refer to \cite{Mollin} for basic aspects of algebra and number theory (including notation and terms not defined here).

\section{Main results}\label{sec:mainresult}

\begin{theorem}\label{th:main0}  
Let $R$ be a ring, $G=\{x_1,\ldots,x_n\}$ a finite subgroup of its group of units, and $f: G^n \to R$ a symmetric homogeneous function of degree $d$ such that there exists $g \in G$ for which $g^d-1$ is regular. Then $f(x_1,\ldots,x_n)=0$.
\end{theorem}
\begin{proof}
Since $G$ is a group, it is then clear that the function $G\to G:x\mapsto gx$ is bijective, hence $\{gx_1,\ldots,gx_n\}=G$. Therefore, using that $f$ is symmetric and homogeneous of degree $d$, we have
\begin{displaymath}
\begin{split}
0=g^{d}f(x_1,\ldots,x_n)-g^{d}f(x_1,\ldots,x_n)&=g^{d}f(x_1,\ldots,x_n)-f(gx_1,\ldots,gx_n) \\
&=g^{d}f(x_1,\ldots,x_n)-f(x_1,\ldots,x_n) \\
&=(g^d-1)f(x_1,\ldots,x_n).
\end{split}
\end{displaymath}
The claim follows from the hypothesis that $g^{d}-1$ is regular.
\end{proof}

We state the following corollary, mainly for future references:
\begin{corollary}\label{cor:casopAofth1}
Let $R$ be a commutative ring, $G$ a finite subgroup of its group of units, and $A=\{a_1,\ldots,a_k\}$ a non-empty multiset of integers for which there exists $g \in G$ such that $g^{s(A)}-1$ is regular. Then $p(A)=0$.
\end{corollary}
\begin{proof}
Let $\Phi: G^n\to R$ be the map defined by
$$
(g_1,\ldots,g_n) \mapsto {\sum_{\substack{x_1,\ldots,x_k \,\in\, \{g_1,\ldots,g_n\} \\ x_1,\ldots,x_k \hspace{1mm}\mathrm{pairwise}\hspace{1mm}\mathrm{distinct}}}}x_1^{a_1}\cdots x_k^{a_k},
$$
where $\Phi(g_1,\ldots, g_n)=0$ whenever there are no pairwise distinct elements $x_1\ldots,x_k$ in the multiset $\{g_1,\ldots,g_n\}$, cf. \S{}~\ref{sec:notation}. Then $\Phi$ is a symmetric homogeneous function of degree $s(A)$. At this point, the claim follows from Theorem~\ref{th:main0}, since there exists $g \in G$ such that $g^{s(A)}-1$ is regular.
\end{proof}

In particular, under the assumptions of Corollary~\ref{cor:casopAofth1}, $\exponent$ does not divide $s(A)$ (indeed, in the opposite case, we would have $g^{s(A)}-1=0$ for all $g \in G$).

It is also easily seen, as a consequence of Theorem~\ref{th:main0}, that: 
\begin{corollary}\label{oldth0}
Let $R$ be an integral domain, $G$ a finite subgroup of its group of units, and $f: G^n \to R$ a symmetric homogeneous function of degree $d$ such that $\exponent$ does not divide $d$. Then $f(x_1,\ldots,x_n)=0$.
\end{corollary}
\begin{proof}
It follows by Theorem~\ref{th:main0} that it is sufficient to show that there exists $g \in G$ for which $g^{d}-1$ is regular. Since $\exponent$ does not divide $d$, there exists $g \in G$ such that $g^{d} \neq 1$, i.e., $g^{d}-1 \neq 0$. The claim follows from the fact that $0$ is the unique non-regular element in $R$.
\end{proof}

This provides a generalization of Pierce's result \cite[Theorem~1]{Pierce}, which corresponds to the case where $R$ is the field $\mathbf{Z}_p$, for some odd prime $p$, and $G$ is the cyclic subgroup of non-zero $n$-th residues modulo $p$. 

In particular, if $R$ is an integral domain and $\exponent$ does not divide $s(A)$, then $p(A)=0$ (we avoid further details, cf. Corollary~\ref{cor:casopAofth1}).

On the other hand, Corollary~\ref{oldth0} does not say anything related to the case where the exponent of $G$ divides $s(A)$. In this regard, we state our main result:
\begin{theorem}\label{th:main1}  
Let $R$ be a commutative ring, $G$ a finite subgroup of its group of units, and $A=\{a_1,\ldots,a_k\}$ a non-empty multiset of integers such that, for every $B\subseteq A$ for which $\lambda$ does not divide $s(B)$, there exists $g \in G$ such that $g^{s(B)}-1$ is regular. Then
\begin{equation}\label{eq:mainclaim}
{\sum_{\substack{x_1,\ldots,x_k \,\in\, G \\ x_1,\ldots,x_k \hspace{1mm}\mathrm{pairwise}\hspace{1mm}\mathrm{distinct}}}}x_1^{a_1}\cdots x_k^{a_k}=\sum_{\mathcal{P} \in \mathscr{P}(A)}\prod_{P \in \mathcal{P}}\bigchi(P).
\end{equation}
\end{theorem}

Note that the result simplifies if $R$ is actually a field:
\begin{corollary}\label{cor:field}  
Let $R$ be a field, $G$ a finite subgroup of its group of units, and $A=\{a_1,\ldots,a_k\}$ a non-empty multiset of integers. Then the identity \eqref{eq:mainclaim} holds.
\end{corollary}

Moreover, hypotheses of Theorem~\ref{th:main1} are verified in the following case:
\begin{corollary}\label{cor:main1}  
Let $R$ be a commutative ring with finitely many non-regular elements, $G$ a finite subgroup of its group of units, and $A=\{a_1,\ldots,a_k\}$ a non-empty multiset of integers such that $\lambda$ divides $s(A)$ and, for each $B\subseteq A$ for which $\lambda$ does not divide $s(B)$, there exists $g \in G$ such that 
\begin{equation}\label{eq:oldminimax}
\frac{\mathrm{ord}(g)}{\mathrm{gcd}(s(B),\mathrm{ord}(g))}\ge |D|+1.
\end{equation}
Then the identity \eqref{eq:mainclaim} holds.
\end{corollary}
However, it follows from \cite[Theorem~I]{Ganesan} that any commutative ring having finitely many non-regular elements and which is not an integral domain is necessarily finite.

\begin{corollary}\label{cor:sqaurefreegeneral}  
Let $R$ be a commutative ring, $G$ a finite subgroup of its group of units, and $A=\{a_1,\ldots,a_k\}$ a non-empty multiset of integers. Let also $P_1,\ldots,P_m$ be prime ideals of $R$ and suppose that: 
\begin{enumerate}[label={\rm (\textsc{c}\arabic{*})}]
\item\label{item:C1} $P_1 \cap \cdots \cap P_m=\{0\}$;
\item\label{item:C2} $P_i+P_j=R$ for all $1\le i<j\le m$;
\item\label{item:C2bis} $|G|=|G/P_1|\cdots |G/P_m|$;
\item\label{item:C3} the exponent of $G/P_i$ is equal to $\lambda$ for all $i=1,\ldots,m$.
\end{enumerate}
Then the identity \eqref{eq:mainclaim} holds.
\end{corollary}

The proof of Theorem \ref{th:main1} follow in \S{} \ref{sec:mainproof}, while Corollaries \ref{cor:field}  -- \ref{cor:sqaurefreegeneral} are proved in \S{} \ref{sec:corollaries}. Some applications and concluding remarks follow in \S{} \ref{sec:concludingrmk}.

\section{Proof of Theorem \ref{th:main1}}\label{sec:mainproof} 

The core of the proof of Theorem \ref{th:main1} is to reduce the problem from sums over distinct entries to sums over single entries. This will be achieved with the aid of convolution and inversion formula in partially ordered sets. Then, sums over free entries will be factorized and reduced sums over a single entry, which will be worked out with the help of Theorem \ref{th:main0}.

\begin{proof}[Proof of Theorem \ref{th:main1}] 
Let $\partition$ be the collection of partitions of $\{1,\ldots,k\}$ partially ordered by refinement, that is, $\mathcal{P} \le \mathcal{Q}$ for some $\mathcal{P},\mathcal{Q} \in \partition$ if and only if for each $P \in \mathcal{P}$ there exists $Q \in \mathcal{Q}$ such that $P\subseteq Q$. We denote by $\mathbf{0}$ and $\mathbf{1}$ its minimum and maximum element, i.e.,
\begin{displaymath}
\mathbf{0}=\{\{1\},\ldots,\{k\}\} \,\,\,\,\text{ and }\,\,\,\,\mathbf{1}=\{1,\ldots,k\}.
\end{displaymath}
We refer to \cite[\S{} 3.10]{Stanley} for a thorough account of basic properties of $\partition$.

Let $\interval$ represent the set of pairs $(\mathcal{P},\mathcal{Q}) \in \partition \times \partition$ such that $\mathcal{P}\le \mathcal{Q}$. In other words, the order interval $[\mathcal{P},\mathcal{Q}]$ is non-empty if and only if $(\mathcal{P},\mathcal{Q}) \in \interval$. Lastly, let $\zeta$ denote the indicator function of $\interval$, that is, $\zeta(\mathcal{P},\mathcal{Q})=1$ whenever $\mathcal{P}\le \mathcal{Q}$, otherwise $\zeta(\mathcal{P},\mathcal{Q})=0$.

At this point, let $\mathscr{F}$ be the set of functions $\interval \to R$, equipped with the convolution product $*$ defined by
$$
(f*g)(\mathcal{P},\mathcal{Q})=\sum_{\mathcal{P} \le \mathcal{R} \le \mathcal{Q}}f(\mathcal{P},\mathcal{R})\,g(\mathcal{R},\mathcal{Q}).
$$
for all $f,g \in \mathscr{F}$ and $\mathcal{P},\mathcal{Q} \in \partition$ with $\mathcal{P} \le \mathcal{Q}$ (in particular, the sum is non-empty).

For each $x=(x_1,\ldots,x_k) \in G^k$, we write $\mathcal{P}_x$ for the partition induced by the equivalence relation $\sim$ on $\{1,\ldots,k\}$ for which $i\sim j$ if and only if $x_i=x_j$. Accordingly, define the functions $\alpha, \beta \in \mathscr{F}$ such that
$$
\alpha(\mathcal{P},\mathcal{Q})=\sum_{x \in G^k, \,\mathcal{P} \le \mathcal{P}_x \le \mathcal{Q}}x_1^{a_1}\cdots x_k^{a_k}
$$
and
$$
\beta(\mathcal{P},\mathcal{Q})=\sum_{x \in G^k, \,\mathcal{P}_x=\mathcal{P}}x_1^{a_1}\cdots x_k^{a_k}
$$
for each $(\mathcal{P},\mathcal{Q}) \in \interval$. In this respect, note that $p(A)=\alpha(\mathbf{0},\mathbf{0})$ and $\alpha(\mathcal{P},\mathcal{P})=\beta(\mathcal{P},\mathcal{Q})$ for all $\mathcal{P},\mathcal{Q} \in \partition$ with $\mathcal{P} \le \mathcal{Q}$.

\begin{claim}\label{claim:convolution}
$\zeta * \beta = \alpha$.
\end{claim}
\begin{proof}
It is enough to observe that, for all $\mathcal{P},\mathcal{Q} \in \partition$ with $\mathcal{P} \le \mathcal{Q}$, it holds
\begin{displaymath}
\begin{split}
(\zeta*\beta)(\mathcal{P},\mathcal{Q})&=\sum_{\mathcal{P} \le \mathcal{R} \le \mathcal{Q}}\beta(\mathcal{R},\mathcal{Q})\\ 
&=\sum_{\mathcal{P} \le \mathcal{R} \le \mathcal{Q}}\,\,\sum_{x \in G^k, \,\mathcal{P}_x=\mathcal{R}}x_1^{a_1}\cdots x_k^{a_k}\\
&=\sum_{x \in G^k, \,\mathcal{P} \le \mathcal{P}_x \le \mathcal{Q}}x_1^{a_1}\cdots x_k^{a_k}=\alpha(\mathcal{P},\mathcal{Q}).
\end{split}
\end{displaymath}
\end{proof}

Let $\mu$ denote the M\H{o}bius function (in $\mathscr{F}$), that is, the inverse of $\zeta$ with respect to the convolution $*$. (The existence of $\mu$ follows by \cite[Proposition 3.6.2]{Stanley}, which however deals with field-valued functions. On the other hand, the proof of the mentioned result relies only on the invertibility of $\zeta(\mathcal{P}, \mathcal{P}) = 1$. Hence, $\mu$ exists also if we consider ring-valued functions.) 

\begin{claim}\label{claim:pAmobius}
$p(A)=\sum_{\mathcal{P} \in \partition}\mu(\mathbf{0},\mathcal{P})\,\alpha(\mathcal{P},\mathbf{1})$.
\end{claim}
\begin{proof}
It follows from Claim \ref{claim:convolution} and the associativity of the convolution $*$ that
\begin{displaymath}
\mu*\alpha=\mu*(\zeta*\beta)=(\mu*\zeta)*\beta=\beta.
\end{displaymath}
In particular, we obtain
\begin{displaymath}
\begin{split}
p(A)=\alpha(\mathbf{0},\mathbf{0})=\beta(\mathbf{0},\mathbf{1})=\sum_{x \in G^k,\,\mathbf{0}\le \mathcal{P} \le \mathbf{1}}\mu(\mathbf{0},\mathcal{P})\,\alpha(\mathcal{P},\mathbf{1}),
\end{split}
\end{displaymath}
which is equivalent to the required identity.
\end{proof}

At this point, it follows by \cite[Example 3.10.4]{Stanley} that
\begin{equation}\label{eq:firstfactor}
\mu(\mathbf{0},\mathcal{P})=\prod_{P \in \mathcal{P}}(-1)^{|P|-1}(|P|-1)!.
\end{equation}

Moreover, for each $\mathcal{P} \in \partition$, define 
$$
\mathcal{P}_A:=\{\{a_i: i \in P\}: P \in \mathcal{P}\}.
$$
Then, we can show that:
\begin{claim}\label{claim:secondfactor}
$\alpha(\mathcal{P},\mathbf{1})=|G|^{|\mathcal{P}|}$ if $\mathcal{P}_A \in \mathscr{P}(A)$, and $\alpha(\mathcal{P},\mathbf{1})=0$ otherwise.
\end{claim}
\begin{proof}
Let us say that $\mathcal{P}=\{P_1,\ldots,P_s\}$. In addition, for each $i=1,\ldots,s$, denote by $y_i$ the common value of the $x_j$s for which $j \in P_i$ and define $b_i:=\sum_{j \in P_i}a_j$ (hence $\{b_1,\ldots,b_s\}=\{s(B):B \in \mathcal{P}_A\}$. Grouping together these $x_j$s, it follows that
\begin{displaymath}
\begin{split}
\alpha(\mathcal{P},\mathbf{1})=\sum_{x \in G^k, \,\mathcal{P} \le \mathcal{P}_x}x_1^{a_1}\cdots x_k^{a_k} 
=\sum_{y \in G^s}y_1^{b_1}\cdots y_s^{b_s}=\prod_{i=1}^s \, \sum_{g \in G}g^{b_i}.
\end{split}
\end{displaymath}

Suppose that $\lambda$ does not divide $b_i$ for some $i=1,\ldots,s$. By hypothesis, there exists $g \in G$ such that $g^{b_i}-1$ is regular, and it follows by Theorem \ref{th:main0} that $\sum_{g \in G}g^{b_i}=0$, hence $\alpha(\mathcal{P},\mathbf{1})=0$.

Otherwise, $\lambda$ divides each $b_i$, that is, $\mathcal{P}_A \in \mathscr{P}(A)$. Since $g^{b_i}=1$ for each $i=1,\ldots,s$ and $g\in G$, we conclude that
$\alpha(\mathcal{P},\mathbf{1})=\prod_{i=1}^s  \sum_{g \in G}1=|G|^s$.
\end{proof}

It follows from Claim \ref{claim:pAmobius}, \eqref{eq:firstfactor}, and Claim \ref{claim:secondfactor}, respectively, that
\begin{displaymath}
\begin{split}
p(A)&=\sum_{\mathcal{P} \in \partition}\mu(\mathbf{0},\mathcal{P})\,\alpha(\mathcal{P},\mathbf{1})\\
&=\sum_{\mathcal{P} \in \partition}	\left(\prod_{P \in \mathcal{P}}(-1)^{|P|-1}(|P|-1)!\right)\,\alpha(\mathcal{P},\mathbf{1})\\
&=\sum_{\mathcal{P} \in \partition,\, \mathcal{P}_A \in \mathscr{P}(A)}	\left(\prod_{P \in \mathcal{P}}(-1)^{|P|-1}(|P|-1)!\right)|G|^{|\mathcal{P}|}. 
\end{split}
\end{displaymath}
By the fact that $|\mathcal{P}|=|\mathcal{P}_A|$ for each $P \in \partition$, we conclude that
$$
p(A)=\sum_{\mathcal{P} \in \partition,\, \mathcal{P}_A \in \mathscr{P}(A)}	\prod_{P \in \mathcal{P}}|G|(-1)^{|P|-1}(|P|-1)!=\sum_{\mathcal{P} \in \mathscr{P}(A)}	\prod_{P \in \mathcal{P}}\,\bigchi(P).
$$
\end{proof}

\section{Proof of Corollaries}\label{sec:corollaries}

\begin{proof}[Proof of Corollary \ref{cor:field}]
According to Theorem \ref{th:main1}, it is enough to verify that, for each subset $B\subseteq A$ such that $\exponent$ does not divide $s(B)$, there exists $g \in G$ such that $g^{s(B)}-1$ is regular. Note that the subgroup of units $g^{s(B)}$, with $g \in G$, contains at least two elements. Since $0$ is the unique non-regular element in $R$, it follows that there exists $g \in G$ such that $g^{s(B)}-1$ is invertible, hence regular.
\end{proof}

\begin{proof}[Proof of Corollary \ref{cor:main1}]
Again, it is enough to check the hypotheses of Theorem \ref{th:main1} hold true. To this aim, note that the subgroup of units $g^{s(B)}$, with $g \in G$, contains 
$$
\max_{g \in G} \, \frac{\mathrm{ord}(g)}{\mathrm{gcd}(s(B),\mathrm{ord}(g))}
$$
distinct elements. In turn, according to \eqref{eq:oldminimax}, this is strictly greater than the number of non-regular elements in $R$. It follows that, for each subset $B\subseteq A$ such that $\exponent$ does not divide $s(B)$, there exists $g \in G$ such that $g^{s(B)}-1$ is regular.
\end{proof}

\begin{proof}[Proof of Corollary \ref{cor:sqaurefreegeneral}]
Define $G_i:=G/P_i$ for each $i=1,\ldots,m$, denote by $\lambda_i$ the exponent of each $G_i$, and note that, since the $P_i$s are prime ideals, each factor ring $R/P_i$ is an integral domain (hence, $G_i$ stands for the projection of $G$ in $R/P_i$). According to Chinese Remainder theorem (and using hypotheses \ref{item:C1} -- \ref{item:C2bis}), we obtain the (surjective) isomorphism
\begin{equation}\label{eq:chineseR}
R\simeq \frac{R}{P_1 \cap \cdots \cap P_m} \simeq \frac{R}{P_1} \times \cdots \times \frac{R}{P_m}.
\end{equation}

At this point, fix a subset $B\subseteq A$ such that $\lambda$ does not divide $s(B)$. Note that, according to Theorem \ref{th:main1}, it would be enough to verify that there exists $g \in G$ such that $g^{s(B)}-1$ is regular. 

Using \eqref{eq:chineseR}, we have $G \simeq G_1 \times \cdots \times G_m$. Hence, denoting by $(g_1,\ldots,g_m) \in G_1\times \cdots \times G_m$ the isomorphic element of $g \in G$, we have equivalently to prove that there exist $g_1 \in G_1,\ldots,g_m \in G_m$ such that $g_i^{s(B)}-1$ is regular in $G_i$ for each $i=1,\ldots,m$. 

Thanks to \ref{item:C3}, each $\lambda_i$ does not divide $s(B)$. In particular, the subgroup $\{g_i^{s(B)}: g_i \in G_i\}$ is not a singleton. To conclude, it is enough to observe that $0$ is the unique non-regular element in the integral domain $R/P_i$, therefore there exists $g_i \in G_i$ such that $g_i^{s(B)}-1$ is regular.
\end{proof}

\section{Applications and Concluding Remarks}\label{sec:concludingrmk}

In this section, we provide some concrete applications of our previous results. At first, we obtain a result of Pierce \cite{Pierce}:
\begin{corollary}\label{cor:pierce}
Let $p$ be an odd prime and $n,k \in \mathbf{N}^+$ such that there are exactly $2k$ non-zero $n$-th residues modulo $p$. Then
\begin{equation}\label{eq:pierceic}
{\sum_{{\{x_1,\ldots,x_k\} \text{ pairwise distinct }n\text{-th residues} \bmod p}}}x_1^{2}\cdots x_k^{2}\equiv 2(-1)^{k-1}\hspace{2mm}\bmod{p}.
\end{equation}
\end{corollary}
\begin{proof}
Set $R=\mathbf{Z}_p$ and $G$ equals to its subgroup of non-zero $n$-th residues. Note that, by hypothesis, it holds
$$
|G|=\lambda=\frac{p-1}{\mathrm{gcd}(n,p-1)}=2k.
$$
In addition, by the fact the permutations of a tuple $(x_1,\ldots,x_k)$ are not counted in \eqref{eq:pierceic}, we have
$$
{\sum_{\{x_1,\ldots,x_k\} \text{ pairwise distinct }n\text{-th residues} \bmod p}}x_1^{2}\cdots x_k^{2}=\frac{1}{k!}\,p(A)
$$
Here, $A$ is the multiset $\{2,\ldots,2\}$, where the $2$ repeats $k$ times. The claim follows by Corollary \ref{cor:field}, indeed $\mathscr{P}(A)=\{\{A\}\}$ so that 
\begin{equation}\label{eq:piercelast}
p(A)\equiv \bigchi(A)= |G|(-1)^{k-1}(k-1)!= 2(-1)^{k-1}k!\,\bmod{p}.
\end{equation}
\end{proof}
However, the above proof reveals that congruence \eqref{eq:piercelast} holds even if that the multiset of exponents $\{2,\ldots,2\}$ is replaced by a multiset of positive integers $A=\{a_1,\ldots,a_k\}$ such that $s(A)=2k$.

We conclude with the following two applications of the ring $\mathbf{Z}_m$ with $m$ prime power and $m$ squarefree, respectively.

\begin{corollary}\label{cor:primepowers}
Let $p^s$ be a power of an odd prime, fix $n\in \mathbf{N}^+$ coprime with $p-1$, and fix $a_1,\ldots,a_k \in \mathbf{Z}$ such that $\varphi(p^s)$ divides $a_1+\cdots+a_k$ and $p-1$ does not divide $\sum_{i \in I}a_i$ for all non-empty proper subsets $I\subseteq \{1,\ldots,k\}$. Then
$$
\sum x_1^{a_1}\cdots x_k^{a_k} \equiv \frac{\varphi(p^s)}{\mathrm{gcd}(n,\varphi(p^s))}(-1)^{k-1}(k-1)! \hspace{1mm}\bmod{p^s},
$$
where the summation is taken over all pairwise distinct $n$-th residues $x_1,\ldots,x_k$ modulo $p^s$ coprime with $p$.
\end{corollary}
\begin{proof}
Set $R=\mathbf{Z}_{p^s}$ and let $G$ be equal to the cyclic subgroup of $R$ of (non-zero) $n$-th residues coprime with $p$. Let $q$ be a primitive root modulo $p^s$ so that $G$ is the subgroup generated by $q^n$ and
$$
|G|=\exponent=\frac{\varphi(p^s)}{\mathrm{gcd}(n,\varphi(p^s))}=\frac{p^{s-1}}{\mathrm{gcd}(n,p^{s-1})}(p-1).
$$

Therefore, $p-1$ divides $\exponent$ which, in turn, divides $\varphi(p^s)$. It follows by the standing assumptions that $\exponent$ divides $s(B)$, for some $B\subseteq A$, if and only if $B=\emptyset$ or $B=A$. In particular, this implies that $\mathscr{P}(A)=\{\{A\}\}$. 

Lastly, fix a non-empty proper subset $B\subseteq A$. We claim that there exists $g \in G$ such that $g^{s(B)}-1$ is not multiple of $p$, hence regular. Indeed, we obtain in $\mathbf{Z}_{p^s}$ that
$$
\{g^{s(B)}-1:g \in G\}=\{q^{s(B)nr}-1: r=1,\ldots,\exponent\}
$$
Since $n$ is coprime with $p-1$ and $p-1$ does not divide $s(B)$, then $p-1$ does not divide $s(B)n$. It follows that there exists a non-zero $n$-th residues modulo $p^s$ which has not remainder $1$ modulo $p$. The claim follows by Theorem \ref{th:main1}, indeed $p(A)\equiv \bigchi(A)=\lambda (-1)^{k-1}(k-1)!\,\bmod{p^s}$.
\end{proof}

\begin{corollary}\label{cor:squarefree}
Let $p_1,\ldots,p_r$ be pairwise distinct odd primes, define $d:=\mathrm{gcd}(p_1-1,\ldots,p_r-1)$, and fix $a_1,\ldots,a_k \in \mathbf{N}^+$ such that $a_1+\cdots+a_k=d$. Then
$$
\sum x_1^{a_1}\cdots x_k^{a_k} \equiv \frac{\varphi(p_1\cdots p_r)}{d^r}(-1)^{k-1}(k-1)! \hspace{1mm}\bmod{p_1\cdots p_r},
$$
where the summation is taken over all pairwise distinct $x_1,\ldots,x_k$ modulo $p_1\cdots p_r$ coprime with $p_1\cdots p_r$ such that each $x_i$ is a $\frac{p_j-1}{d}$-th residue modulo $p_j$ for all $j=1,\ldots,r$.
\end{corollary}
\begin{proof}
We will verify that the hypotheses of Corollary \ref{cor:sqaurefreegeneral} hold. In this regard, set $R=\mathbf{Z}_{p_1\cdots p_r}$ and $G$ its subgroup of (non-zero) $d$-th residues modulo $p_1\cdots p_r$ coprime with $p_1\cdots p_r$. Moreover, denote by $P_i$ the prime ideal $p_i\mathbf{Z}$ for each $i=1,\ldots,r$. It is straighforward to check the conditions \ref{item:C1} -- \ref{item:C2bis} hold. In addition, we have $\exponent=d$, indeed $x^d \equiv 1\pmod{p_i}$ for each $x \in G$ and $i=1,\ldots,r$. On the other hand, $\lambda_i=d$ for each $i=1,\ldots,r$, hence \ref{item:C3} is also verified.

Lastly, note that $\mathscr{P}(A)=\{\{A\}\}$, which allows us to conclude that
\begin{displaymath}
\begin{split}
p(A) &\equiv \bigchi(A)=|G|(-1)^{k-1}(k-1)!\\
&=|G_1|\cdots |G_r|(-1)^{k-1}(k-1)!=\frac{\varphi(p_1\cdots p_r)}{d^r}(-1)^{k-1}(k-1)! \hspace{1mm}\bmod{p_1\cdots p_r}.
\end{split}
\end{displaymath}
\end{proof}

To conclude the section, define
$$
A^\free:=\{a \in A: \exponent \text{ does}\text{ not}\text{ divide }a\}\,\,\,\text{ and }\,\,\,\ell:=|A^\free|.
$$
\begin{claim}\label{claim:Afree}
If $A^\free \neq \emptyset$ then $p(A)=(n-\ell)(n-\ell-1)\cdots (n-k+1) p(A^\free)$.
\end{claim}
\begin{proof}
It is enough to note that, if $A^\free=\{a_1,\ldots,a_{\ell}\}$ and $\ell \in \mathbf{N}^+$, then each summand $x_1^{a_1}\cdots x_{\ell}^{a_{\ell}}$ in $p(A^\free)$ appears in \eqref{eq:pA} exactly $(n-\ell)!/(n-k)!$ times.
\end{proof}
Accordingly, if the assumptions the Theorem \ref{th:main1} hold and $A^\free\neq \emptyset$, then
\begin{displaymath}
\begin{split}
p(A)&=\frac{(n-\ell)!}{(n-k)!}\sum_{\mathcal{P} \in \mathscr{P}(A^\free)}\prod_{P \in \mathcal{P}}\bigchi(P)\\
&=\frac{(n-\ell)!}{(n-k)!}\sum_{\mathcal{P} \in \mathscr{P}(A^\free)}n^{|\mathcal{P}|}\prod_{P \in \mathcal{P}}(-1)^{|P|-1}(|P|-1)!\\
&=(-1)^k\frac{(n-\ell)!}{(n-k)!}\sum_{\mathcal{P} \in \mathscr{P}(A^\free)}(-n)^{|\mathcal{P}|}\prod_{P \in \mathcal{P}}(|P|-1)!.
\end{split}
\end{displaymath}

\section{Acknowledgements}\label{sec:acknowledgments} 

We are grateful to professor Lo\"ic Greni\'{e} (University of Bergamo, Italy) for a careful proof-reading of the manuscript and to Salvatore Tringali (University of Graz, Austria) and Pierfrancesco Carlucci (Università di Roma Tor Vergata, Italy) for helpful comments. We thank also the anonymous referee for many suggestions which improved greatly the readability of the article.

\end{document}